\documentclass{article}

\usepackage[english]{babel}

\usepackage[letterpaper,top=2cm,bottom=2cm,left=3cm,right=3cm,marginparwidth=1.75cm]{geometry}

\usepackage{amsmath,amsfonts,amssymb, stmaryrd,mathtools,amsthm}
\usepackage{graphicx}
\usepackage{natbib}
\usepackage[colorlinks=true, allcolors=blue]{hyperref}
\allowdisplaybreaks[4]

\newtheorem{theorem}{Theorem}[section]

\newtheorem{corollary}[theorem]{Corollary}
\newtheorem{lemma}[theorem]{Lemma}
\newtheorem{proposition}[theorem]{Proposition}

\newtheorem{remark}[theorem]{Remark}

\newcommand{\n}{\nabla}
\newcommand{\no}{\mathbf{n}}
\newcommand{\D}{\overline{\nabla}}

\newcommand{\R}{\mathbb{R}}

\title{ Harnack Inequality for $f$-Mean Curvature Flow}
\author{ Xiang-Dong Li \footnote{1. State Key Laboratory of Mathematical Sciences, Academy of Mathematics
and Systems Science, Chinese Academy of Sciences, No. 55, Zhongguancun East Road, Beijing,
100190, China; 2. School of Mathematical Sciences, University of Chinese Academy of Sciences,
Beijing, 100049, China. Research supported by National Key R\&D Program of China (No. 2020YF0712702), NSFC No.
 12171458, and Key Laboratory RCSDS, CAS, No. 2008DP173182.  \texttt{xdli@amt.ac.cn}}  \quad and \quad Qi Yan \footnote{Academy of Mathematics
and Systems Science, Chinese Academy of Sciences, No. 55, Zhongguancun East Road, Beijing,
100190.  \texttt{yanqi19@mails.ucas.ac.cn}}     }

\begin{document}
\maketitle

\begin{abstract}
In this paper, we prove a Li-Yau-Hamilton type Harnack estimate for the $f$-mean curvature flow in Euclidean space, which can be viewed as a gradient flow of the weighed area functional with the measure density function $e^{-f}$.
\end{abstract}

\section{Introduction}
\subsection{$f$-Mean Curvature Flow}
Mean curvature flow is an important research object in differential geometry, with significant theoretical and applied meanings. In physics, as the sharp interface limit of the Allen-Cahn equation for binary alloy phase fields, it describes the motion of biphasic interfaces, and can also be used to describe crystal growth processes.

Given a family of $n$-dimensional hypersurfaces $M_t$ in $\mathbb{R}^{n+1}$, it is called mean curvature flow if it satisfies
\[
\partial_t x=-H\mathbf{n}=\mathbf{H}.
\]
where $x$ denotes the position vector in $\mathbb{R}^{n+1}$, $\mathbf{n}$ is the unit outer normal vector at $x\in M_t$, $H$ is the mean curvature of the hypersurface $M_t$ at $x$, and $\mathbf{H}=-H\no$ is the mean curvature vector.

Mean curvature flow is the gradient flow of the area functional (in the sense of the $L^2$ metric). Let $M_t$ be a family of hypersurfaces with $M_0=M$, and consider the area functional $\mathrm{Area}$ as defined on all $n$-dimensional hypersurfaces in $\mathbb{R}^{n+1}$. Then, by the first variation formula, we have
\[
  \frac{d}{d t} \operatorname{Area}\left(M_t\right)=\int_{M_t}\left\langle\partial_t x, H \mathbf{n}\right\rangle.
\]
This implies
\[
  \nabla \mathrm{Area}=H \mathbf{n}.
\]
Therefore,
\[
  \partial_t x=-H\mathbf{n}=-\nabla \mathrm{Area}.
\]
Thus, mean curvature flow can be regarded as the gradient flow of the area functional.

The $f$-mean curvature flow is a generalization of mean curvature flow. From the perspective of gradient flow, we can introduce the so-called $f$-mean curvature flow. We consider a function $f$ defined on $\mathbb{R}^{n+1}$, and a weighted measure $d\mu=e^{-f}dx$, where $dx$ is the Lebesgue measure. The weighted area functional with respect to this weighted measure is defined as
\[
  \mathrm{Area}_f(M_t)=\int_{M_t}e^{-f}dx,
\]
By the first variation formula, we obtain
\[
  \frac{d}{d t} \operatorname{Area}_f\left(M_t\right)=\int_{M_t}\left\langle\partial_t x, H_f \mathbf{n}\right\rangle d\mu,
\]
where $H_f=H-\langle \overline{\nabla} f,\mathbf{n}\rangle$ is called the \textbf{$f$-mean curvature} (see \cite{Wei09,LW2015}), here $\overline{\nabla}$ is the gradient in $\R^{n+1}$. The gradient flow of the weighted area functional
\[
  \partial_t x=-H_f\mathbf{n}
\]
is called the \textbf{$f$-mean curvature flow} (see \cite{JL08}). Under the above setting, the $f$-mean curvature flow can be viewed as the mean curvature flow on the smooth metric measure space $(\R^{n+1},dx^2,d\mu)$. When $f=\frac{C|x|^2}{2}$ (where $C>0$), this curvature flow is called the Gauss mean curvature flow, see \citep{BM10,BM14}.
When $f$ is constant, the $f$-mean curvature flow is just the mean curvature flow.

The $f$-mean curvature flow can also be viewed as a mean curvature flow. Let $F_t:M^n\to\R^{n+1}$ be the $f$-mean curvature flow, and $dx^2$ be the Euclidean metric on $\R^{n+1}$. Consider the following mapping:
\begin{equation*}\label{flow2}
	\begin{array}[pos]{rcl}
		\widetilde{F}_t:M^n\times\R&\to&\R^{n+1}\times\R\\
		(x,s)&\mapsto&(F_t(x),s)
	\end{array}
\end{equation*}
where the metric on $\R^{n+1}\times\R$ is $dx^2+e^{2f(x)}ds^2.$ Let $\widetilde{M}_t=M_t\times\R$. From \citep[Corollary 2.4]{Smoczyk01}, we know that
\[
	\mathbf{H}_{\widetilde{M}_t}(x,s)=\left(\mathbf{H}_f(x),0\right),
\]
where $\mathbf{H}_f=-H_f\no$ is the $f$-mean curvature vector of $M_t=F_t(M)$, and $\mathbf{H}_{\widetilde{M}_t}$ is the mean curvature vector of $\widetilde{M}_t$ in the Riemannian manifold $(\mathbb{R}^{n+1}\times\R,dx^2+e^{2f(x)}ds^2)$. From this, we obtain
\[
	\frac{\partial\widetilde{F}}{\partial t}(x,s)=\left(\frac{\partial F}{\partial t},0\right)=\left(\mathbf{H}_f,0 \right)=\mathbf{H}_{\widetilde{M}_t}(x,s).
\]
Therefore, the $f$-mean curvature flow is the mean curvature flow in the Riemannian manifold $(\mathbb{R}^{n+1}\times\R,dx^2+e^{2f(x)}ds^2)$.

The $f$-mean curvature flow also originates from the study of Ginzburg-Landau vortex flows as follows.
\[
\frac{\partial V_{\epsilon}}{\partial t}=\overline{\Delta}V_{\epsilon}+\overline{\nabla}f\overline{\nabla}V_{\epsilon}+AV_{\epsilon}+\frac{BV_{\epsilon}}{\epsilon^2}(1-|V_{\epsilon}|^2),\quad V_{\epsilon}:\mathbb{R}^n\times [0,T]\to\mathbb{R}^2.
\]
See \cite{JL06}.
In fact, mean curvature flow is closely related to the Allen-Cahn equation. The Allen-Cahn equation is a real-valued simplification of the Ginzburg-Landau theory, describing the dynamics of bistable systems (such as phase separation)
\[
\partial_t u=\Delta u-\frac{1}{\epsilon^2}W^{\prime}(u),\quad u:\mathbb{R}^n\times [0,T]\to\mathbb{R}.
\]
The level set $M_t:=\{x\in\mathbb{R}^n: u(x,t)=0\}$ of the limit $u(t)$ of the solution $u_{\epsilon}(t)$ of the Allen-Cahn equation as $\epsilon\to 0$ is a mean curvature flow. The level set of the limit $V(t)$ of the solution $V_{\epsilon}(t)$ of the Ginzburg-Landau equation as $\epsilon\to 0$ is a codimension-2 $f$-mean curvature flow.

R. Hamilton, in his 1995 paper \citep{Hamilton95}, obtained a Li-Yau type Harnack inequality for the mean curvature under mean curvature flow in Euclidean space:
\begin{theorem}[Hamilton 95']\label{Harnack inequality for mean curvature flow}
Every weakly convex solution of the mean curvature flow
\begin{equation}
    \partial_tF=-H\no
\end{equation}
satisfies: for $t>0$, and for any tangent vector $V$, the following inequality holds:
\begin{equation}\label{mean curvature flow Harnack inequality}
    \partial_t H+2DH(V)+h(V,V)+\frac{H}{2t}\geq 0.
\end{equation}
where $h$ is the second fundamental form.
\end{theorem}
In fact, Ricci flow as well as Gaussian curvature flow have similar Harnack estimate property, see \cite{Hamilton95}and \cite{Chow91} respectively. Knot Smoczyk  proved in \cite{Smoczyk97} Harnack inequality for parabolic flows of compact orientable  hypersurface in $\R^{n+1}$, where the normal velocity is given by a smooth function depending only on the mean curvature .  Bryan and Ivaki proved in \cite{BI20} Harnack inequality for mean curvature flow on the sphere.

In this paper, we will investigate the Harnack inequality for the $f$-mean curvature flow.

\section{Evolution Equations for the $f$-Mean Curvature Flow}
In this chapter, we always use the following notation: $\langle\cdot,\cdot\rangle$ denotes the standard inner product in $\mathbb{R}^{n+1}$, and let $F: M^n\to\mathbb{R}^{n+1}$ be the $f$-mean curvature flow:
\begin{equation}\label{f-mean curvature flow}
    \dfrac{\partial F}{\partial t}=-H_fN.
\end{equation}
The induced metric and the second fundamental form on $M$ are defined as
\[
	g_{ij}(x)=\left\langle\frac{\partial F}{\partial x_i},\frac{\partial F}{\partial x_j}\right\rangle(x),\quad 
	h_{ij}(x)=-\left\langle N,\dfrac{\partial^2 F}{\partial x_i\partial x_j}\right\rangle(x),\quad x\in M.
\] 
Using the Einstein summation convention, the traces related to the second fundamental form are as follows:
\begin{gather*}
	H=g^{ij}h_{ij}=\operatorname{tr}h,\quad 
	|h|^2=g^{ij}g^{kl}h_{ik}h_{jl},\\
	C=g^{ij}g^{kl}g^{mn}h_{ik}h_{lm}h_{nj}=\operatorname{tr}h^3,\quad 
	Z=HC-|h|^4.
\end{gather*}
Additionally, the Gauss and Codazzi equations are:
\begin{gather}\label{Gauss-Codazzi}
	R_{ijkl}=h_{ik}h_{jl}-h_{il}h_{jk}, \\
	 \nabla_i h_{jk}= \nabla_{j}h_{ik}.
\end{gather}

For Euclidean hypersurfaces, we have the following basic properties:
\begin{proposition}
	For any hypersurface $F:M\to\mathbb{R}^{n+1}$, we have
	\begin{equation}
		\begin{aligned}
		\nabla_i \nabla_j F & =-h_{i j} N, \\
		\nabla_i N & =h_i^l \nabla_l F, \\
		\nabla_i \nabla_j H & =\Delta h_{i j}-H h_i^l h_{l j}+|h|^2 h_{i j}, \\
		2 h^{i j} \nabla_i \nabla_j H & =\Delta|h|^2-2|\nabla h|^2-2 Z,
		\end{aligned}
	\end{equation}
	where $\nabla$ denotes the Riemannian connection on $M$, and $\Delta$ is the Laplace-Beltrami operator on $M$.
\end{proposition}
The proof of the following proposition can be found in \cite{JL08}.
\begin{proposition}
	Suppose the $f$-mean curvature flow \eqref{f-mean curvature flow} holds on $[0,T)$, $T<\infty$. Then
	\begin{flalign}
		\dfrac{\partial g_{ij}}{\partial t}&=-2H_fh_{ij}, \\
		\dfrac{\partial N}{\partial t}&=\nabla H_f,\\
		\dfrac{\partial H_f}{\partial t}&=\Delta H_f-\langle\overline{\nabla}f,\nabla H_f\rangle+\left(|h|^2+\overline{\nabla}^2f(N,N)\right)H_f,\\
		\dfrac{\partial h_{ij}}{\partial t}&=\nabla_i\nabla_j H_f-H_fh_i^lh_{lj}.
	\end{flalign}
\end{proposition}
Using the maximum principle, we can obtain
\begin{proposition}
    For the $f$-mean curvature flow \eqref{f-mean curvature flow}, if $H_f\geq 0$ at $t=0$, then $H_f\geq 0$ for all $t\geq 0$.
\end{proposition}

H.Y. Jian and Y. Liu, in the literature \cite{JL08}, gave more precise evolution equations for $h_{ij}$ and used Hamilton's tensor maximum principle to prove that when $h_{ij}(0)\geq 0$ and $f$ satisfies $\overline{\n}^3 f\equiv 0$, $h_{ij}$ remains positive definite.
\begin{proposition}
	Suppose the $f$-mean curvature flow \eqref{f-mean curvature flow} holds on $[0,T)$, where $T<\infty$, then we have
	\begin{equation}
		\begin{aligned}
			\dfrac{\partial h_{ij}}{\partial t}=
			& \Delta h_{ij}-\langle\overline{\nabla}f,\nabla h_{ij}\rangle-2H_fh^l_ih_{lj}+\left(|h|^2+\overline{\nabla}^2f(N,N)\right)h_{ij}\\
			&-h^l_i\overline{\nabla}^2f(\nabla_jF,\nabla_lF)-h^l_j\overline{\nabla}^2f(\nabla_iF,\nabla_lF)-\overline{\nabla}^3f(N,\nabla_i F,\nabla_jF).
		\end{aligned}
	\end{equation}
	
	If $f$ satisfies $\overline{\n}^3f\equiv 0$\footnote{We attempted to weaken this condition to: $\overline{\n}^3f(U,V,V)\leq 0$ for all vectors $U,V\in\R^{n+1}$ satisfying $\langle U, V\rangle=0$, but this is equivalent to $\overline{\n}^3f\equiv 0$.}, then if $h_{ij}$ is nonnegative at $t=0$, it remains nonnegative for all $t\geq 0$.
\end{proposition}

\section{Differential Harnack Inequality for the $f$-Mean Curvature Flow}
This section will prove the Harnack inequality for the $f$-mean curvature flow.

\subsection{Notation System}
To simplify calculations, we use moving orthonormal frames.

Let $U\subset M^n$ be a coordinate neighborhood with local coordinates $\{x^i\}$, and $Y=\{Y_a\}$ be an orthonormal frame, where $Y_1,\dots,Y_n$ are tangent vectors to $M$, and $Y_a=y^i_a\frac{\partial }{\partial x^i}$. This defines coordinates $\{x^i,y^j_a\}$ on the frame bundle, so $\{\frac{\partial }{\partial x^i},\frac{\partial}{\partial y^j_a}\}$ form a basis for the bundle tangent vectors.

Let $\n_a$ be the horizontal vector field on the frame bundle, projecting to $Y_a$ at the point $Y$; $\n^a_b$ be the vertical vector field on the frame bundle fibers, representing the action of the general linear group $\mathrm{GL}(n)$. Through these vector fields, any tensor induces a system of functions on the frame bundle. For example:
For a covector $V:TM\to\mathbb{R}$:
\[
	V=\{V_a\}\quad\text{where}\quad V(Y_a)=V_a.
\]
For a tensor $V:TM\times TM\to TM$:
\[
	V=\{V^c_{ab}\}\quad \text{where}\quad V(Y_a,Y_b)=V^c_{ab}Y_c.
\]
Covariant derivatives can be understood as $\n_a$ acting on the tensor component functions:
When $V=\{V_a\}$:
\[
	\n V=\{\n_aV_b\}\quad\text{where}\quad \n V(Y_a)(Y_b)=\n _aV_b.
\]
When $V=\{V^d_{bc}\}$:
\[
	\n V=\{\n_aV^d_{bc}\}\quad\text{where}\quad \n V(Y_a)(Y_b,Y_c)=\n_aV^d_{bc}Y_d.
\]
In local coordinates $\{x^i, y^j_a\}$, the frame bundle vector fields are:
\[
	\n_a=y^i_a\left(\dfrac{\partial}{\partial x^i}-\Gamma^k_{ij}y^j_b\dfrac{\partial }{\partial y^k_b}\right), \quad \n^a_b=y^i_b\dfrac{\partial}{\partial y^i_a}.
\]

It is easy to compute the action of $\n^a_b$ on tensors:
For a function $\phi$ on $M$, $\n^a_b\phi=0.$ For a 1-form $V$, $\n^a_bV_c=\delta^a_cV_b.$
For a 2-tensor $W$, $\n^a_bW_{cd}=\delta^a_cW_{bd}+\delta^a_dW_{cb}.$
  
We have the following commutation relations:
\begin{equation}
	\begin{aligned}
		\n_a\n_b-\n_b\n_a&=R^c_{abd}\n^d_c.\\
		\n^a_b\n_c-\n_c\n^a_b&=\delta^a_c\n_b.
	\end{aligned}
\end{equation}

In calculations, we use the orthonormal frame bundle, where the Riemannian metric $g_{ab}=g(Y_a,Y_b)=\delta_{ab}$. Since the orthonormal frame depends on time, for the $f$-mean curvature flow \eqref{f-mean curvature flow}, we define the covariant time derivative:
\[
	\n_t=\partial_t+H_fh_{ab}g^{bc}\n^a_c=\partial_t+H_fh_{ac}\n^a_c.
\]
Then $\n_tg_{ab}=0$. One can directly verify:
\begin{equation}\label{nabla t,a}
	\begin{aligned}
		\relax[\n_t,\n_a]&=H_fh_{ab}\n_b+\left(h_{ac}\n_bH_f-h_{ab}\n_cH_f\right)\n^b_c.\\
		[\n_t,\n^a_b]&=0.
	\end{aligned}
\end{equation}

To simplify notation, define:
\[
	\n^2_{ab}f=\D^2f(\n_aF,\n_bF), \quad \D^3_{ab}f(N)=\D^3f(N,\n_aF,\n_bF).
\]
Define the differential operator on tensors:
\begin{equation}
	L=\Delta-\langle\D f,\n\rangle-\n^2_{ab}f\n^a_b.
\end{equation}
For any function $\phi$, $L\phi=\Delta\phi-\langle\D f,\n\phi\rangle$.

Since $\n_af =\langle\D f,\n_aF\rangle$ and $\n=\n_aF\n_a$, we can rewrite as:
\begin{equation}
	L:=\Delta-\n_a f\n_a-\n^2_{ab}f\n^a_b.
\end{equation}
Using the above notation, the evolution equations for $H_f$ and $h_{ab}$ can be expressed as:
\begin{flalign}
	\n_tH_f&=L H_f+\left(|h|^2+\D^2f(N,N)\right)H_f.\\
	\n_th_{ab}&=\n_a\n_bH_f+H_fh^2_{ab}=L h_{ab}+\left(|h|^2+\D^2f(N,N)\right)h_{ab}-\D^3_{ab}f(N). 
\end{flalign}

\subsection{Three Key Commutation Relations}
For a smooth function $\phi\in C^{\infty}(M)$, we have the following commutation relations:
\begin{lemma}
	\begin{equation}
		\left(L\n_a-\n_aL\right)\phi=H_fh_{ab}\n_b\phi-h^2_{ab}\n_b\phi.
	\end{equation}
\end{lemma}
\begin{proof}
	On the hypersurface, we have:
	\[
		(\Delta\n_a-\n_a\Delta)\phi=R_{ab}\n_b\phi=Hh_{ab}\n_b\phi-h^2_{ab}\n_b\phi.
	\]
	By computation:
	\begin{align*}
		\dfrac{\partial}{\partial x^i}\n\phi-\n\dfrac{\partial}{\partial x^i}\phi 
		& =\dfrac{\partial}{\partial x^i}\left(\dfrac{\partial\phi}{\partial x^j}g^{jk}\dfrac{\partial F}{\partial x^k}\right)-\dfrac{\partial^2\phi}{\partial x^j\partial x^i}g^{jk}\dfrac{\partial F}{\partial x^k}\\
		& =\dfrac{\partial\phi}{\partial x^j}g^{jk}\dfrac{\partial^2F}{\partial x^i\partial x^k}\\
		& =-h_{ik}\dfrac{\partial\phi}{\partial x^j}g^{jk}N.
	\end{align*}
	Hence:
	\[
		\langle\D f,\n_a\n\phi-\n\n_a\phi\rangle=-\langle\D f,N\rangle h_{ab}\n_b\phi.
	\]
	Note that:
	\[
		\n^2_{cd}f(\n_a\n^d_c-\n^d_c\n_a)\phi=-\n^2_{ac}f\n_c\phi.
	\]
	And:
	\[
		\langle\n_a(\D f),\n\phi\rangle=\D^2f(\n_aF,\n\phi)=\n^2_{ab}f\n_b\phi.
	\]
	Finally:
	\begin{align*}
		\left(L\n_a-\n_aL\right)\phi
		&=(\Delta\n_a-\n_a\Delta)\phi+\langle\D f,\n_a\n\phi-\n\n_a\phi\rangle+\n^2f_{cd}
			\left(\n_a\n^d_c-\n^d_c\n_a\right)\phi\\
		&\quad+\langle\n_a(\D f),\n\phi\rangle+\D^2f(\n_a\n_cF,\n_dF)\n^d_c\phi+\D^2f(\n_cF,\n_a\n_dF)\n^d_c\phi.\\
		&=H_fh_{ab}\n_b\phi-h^2_{ab}\n_b\phi.
	\end{align*}
\end{proof}
This lemma immediately implies:
\begin{lemma}
	\[
		(\n_t-L)\n_a\phi-\n_a(\n_t-L)\phi=h^2_{ab}\n_b\phi.
	\]
\end{lemma}

\begin{lemma}
	\begin{align}
		(\n_t-L)L\phi-L(\n_t-L)L\phi
		&=2H_fh_{ab}\n_a\n_b \phi+2h_{ab}\n_aH_f\n_b\phi+2H_f\D^2f(\n_aF,N)\n_a\phi.
	\end{align}
\end{lemma}
\begin{proof}
	Apply (\ref{nabla t,a}) to $\n_a\phi$ and $\phi$ respectively:
	\begin{align*}
		\n_t\n_a\n_a\phi-\n_a\n_t\n_a\phi
		&=[\n_t,\n_a]\n_a\phi=H_fh_{ab}\n_b\n_a\phi+(h_{ac}\n_bH_f-h_{ab}\n_cH_f)\n^d_c\n_a\phi\\
		&=H_fh_{ab}\n_a\n_b\phi+(h_{ac}\n_bH_f-h_{ab}\n_cH_f)\delta^b_a\n_c\phi\\
		&=H_fh_{ab}\n_a\n_b\phi+h_{ab}\n_aH_f\n_b\phi-H\n_aH_f\n_a\phi.
	\end{align*}
	And:
	\[
		\n_t\n_a\phi-\n_a\n_t\phi=H_fh_{ab}\n_b\phi.
	\]
	Apply $\n_a$ to both sides:
	\begin{align*}
		\n_a\n_t\n_a\phi-\n_a\n_a\n_t\phi
		&=\n_a(H_fh_{ab}\n_b\phi)=h_{ab}\n_aH_f\n_b\phi+H_f\n_ah_{ab}\n_b\phi+H_fh_{ab}\n_a\n_b\phi\\
		&=h_{ab}\n_aH_f\n_b\phi+H_f\n_aH\n_a\phi+H_fh_{ab}\n_a\n_b\phi.
	\end{align*}
	Therefore:
	\begin{align*}
		\n_t\Delta\phi-\Delta\n_t\phi
		&=\n_t\n_a\n_a\phi-\n_a\n_a\n_t\phi\\
		&=2H_fh_{ab}\n_a\n_b\phi+2h_{ab}\n_aH_f\n_b\phi+H_f\n_aH\n_b\phi-H\n_aH_f\n_b\phi.
	\end{align*}
	Combining with:
	\begin{align*}
		\n_t\n_af\n_a\phi 
		&=\n_a\n_tf\n_a\phi+H_fh_{ab}\n_af\n_b\phi=\n_a\langle\D f,-H_fN\rangle\n_a\phi+H_fh_{ab}\n_af\n_b\phi\\
		&=-H_f\D^2f(\n_aF,N)\n_a\phi-\langle\D f,N\rangle\n_aH_f\n_a\phi-H_fh_{ab}\n_bf\n_a\phi+H_fh_{ab}\n_af\n_b\phi\\
		&=-H_f\D^2f(\n_aF,N)\n_a\phi-\langle\D f,N\rangle\n_aH_f\n_a\phi.
	\end{align*}
	We get:
	\begin{align*}
		(\n_t-L)L\phi-L(\n_t-L)\phi
		&=\n_t\Delta\phi-\Delta\n_t\phi+\n_af(\n_a\n_t\phi-\n_t\n_a\phi)-\n_t\n_af\n_a\phi\\
		&=2H_fh_{ab}\n_a\n_b\phi+2h_{ab}\n_aH_f\n_b\phi+H_f\n_aH\n_b\phi-H\n_aH_f\n_b\phi\\
		&\quad-H_fh_{ab}\n_af\n_b\phi+H_f\D^2f(\n_aF,N)\n_a\phi+\langle\D f,N\rangle\n_aH_f\n_a\phi\\
		&=2H_fh_{ab}\n_a\n_b \phi+2h_{ab}\n_aH_f\n_b\phi+H_f\n_a\langle\D f,N\rangle\n_a\phi\\
		&\quad -H_fh_{ab}\n_af\n_b\phi+H_f\D^2f(\n_aF,N)\n_a\phi\\
		&=2H_fh_{ab}\n_a\n_b \phi+2h_{ab}\n_aH_f\n_b\phi+2H_f\D^2f(\n_aF,N)\n_a\phi.
	\end{align*}
\end{proof}

\subsection{Evolution Equation of the Harnack Quantity}
Define the following tensors, where $V_a$ is any tangent vector on $M$:
\begin{align*}
	X_a&=\n_aH_f+h_{ab}V_b,\\
	Y_{ab}&=\n_aV_b-H_fh_{ab},\\
	W_{ab}&=\n_th_{ab}+V_c\n_ch_{ab},\\
	U_a&=(\n_t-L)V_a+h_{ab}\n_bH_f+2\D^2f(\n_aF,N),\\
	Z&=\n_tH_f+2V_a\n_aH_f+h_{ab}V_aV_b.
\end{align*}
We call $Z$ the \textbf{Harnack quantity} for the $f$-mean curvature flow. Note:
\[
	Z=\partial_tH_f+V_a\n_aH_f+X_aV_a.
\]
We compute the evolution equations for the three terms in $Z$ separately.
First, compute the first term:
\begin{align*}
	(\n_t-&L)(\partial_t H_f)\\
	&=(\n_t-L)\left(LH_f+(|h|^2+\D^2f(N,N))H_f\right)\\
	&=L(\n_t-L)H_f+2H_fh_{ab}\n_a\n_bH_f+2h_{ab}\n_aH_f\n_bH_f+2H_f\D^2f(\n_aF,N)\n_aH_f\\
	&\qquad+(\n_t-L)\left((|h|^2+\D^2f(N,N))H_f\right)\\
	&=\n_t\left((|h|^2+\D^2f(N,N))H_f\right)+2H_fh_{ab}\n_a\n_bH_f+2h_{ab}\n_aH_f\n_bH_f\\
	&\qquad +2H_f\D^2f(\n_aF,N)\n_aH_f\\
	&=\left(|h|^2+\D^2f(N,N)\right)(\partial_tH_f)+2H_fh_{ab}\n_a\n_bH_f+2h_{ab}\n_aH_f\n_bH_f\\
	&\qquad +4H_f\D^2f(\n_aF,N)\n_aH_f+2H_fh_{ab}\n_th_{ab}-H^2_f\D^3f(N,N,N).
\end{align*}
Next, compute the second term:
\begin{align*}
	(\n_t-&L)(V_a\n_aH_f)\\
	&=\n_aH_f(\n_t-L)V_a+V_a(\n_t-L)\n_aH_f-2\n_cV_a\n_c\n_aH_f\\
	&=V_a\n_a\left((|h|^2+\D^2f(N,N))H_f\right)+h^2_{ab}V_a\n_bH_f-2\n_cV_a\n_c\n_aH_f+\n_aH_f(\n_t-L)V_a\\
	&=\left(|h|^2+\D^2f(N,N)\right)(V_a\n_aH_f)+h^2_{ab}V_a\n_bH_f-2\n_cV_a\n_c\n_aH_f+\n_aH_f(\n_t-L)V_a\\
	&\qquad +2H_fh_{bc}V_a\n_ah_{bc}+2V_a\D^2f(\n_bF,N)h_{ab}+H_fV_a\D^3f(\n_aF,N,N).
\end{align*}
Before computing the evolution of the last term, first compute how $X_a$ evolves:
\begin{align*}
	(\n_t-L)X_a
	&=(\n_t-L)\n_aH_f+(\n_a-L)h_{ab}V_b+h_{ab}(\n_t-L)V_b-2\n_ch_{ab}\n_cV_b\\
	&=\n_a\left((|h|^2+\D^2f(N,N))H_f\right)+h^2_{ab}\n_bH_f+\left(|h|^2+\D^2f(N,N)\right)h_{ab}V_b\\
	&\qquad -\D^3f(N,\n_aF,\n_bF)V_b+h_{ab}(\n_t-L)V_b-2\n_ch_{ab}\n_cV_b\\
	&=\left(|h|^2+\D^2f(N,N)\right)X_a+h^2_{ab}\n_bH_f+h_{ab}(\n_t-L)V_b-2\n_ch_{ab}\n_cV_b\\
	&\qquad +2H_fh_{bc}\n_ah_{bc}+2H_f\D^2f(\n_bF,N)h_{ab}+H_f\D^3f(\n_aF,N,N)\\
	&\qquad -\D^3f(N,\n_aF,\n_bF)V_b.
\end{align*}
Now compute the last term:
\begin{align*}
	(\n_t-L)&(X_aV_a)\\
	&=(\n_t-L)X_aV_a+X_a(\n_t-L)V_a-2\n_cV_a\n_cX_a\\
	&=\left(|h|^2+\D^2f(N,N)\right)(X_aV_a)+h^2_{ab}V_a\n_bH_f+h_{ab}V_a(\n_t-L)V_b-2V_a\n_ch_{ab}\n_cV_b\\
	&\qquad +2H_fV_ah_{bc}\n_ah_{bc}+2H_fh_{ab}V_a\D^2f(\n_bF,N)+X_a(\n_t-L)V_a-2\n_cV_a\n_c\n_aH_f\\
	&\qquad -2\n_cV_a\n_ch_{ab}V_b-2\n_cV_ah_{ab}\n_cV_b+H_fV_a\D^3f(\n_aF,N,N)\\
	&\qquad-\D^3f(N,\n_aF,\n_bF)V_aV_b.
\end{align*}

Combining the above three terms, we get the evolution equation for $Z$:
\begin{align*}
	&(\n_t-L)Z\\
	=&\left(|h|^2+\D^2f(N,N)\right)Z\\
	&\qquad+2X_a(\n_t-L)V_a+2(\n_aH_f+h_{ac}V_c)h_{ab}\n_bH_f+4H_f\D^2f(\n_aF,N)X_a\\
	&\qquad+2H_fh_{ab}\left(\n_a\n_bH_f+2V_c\n_ch_{ab}+\n_th_{ab}\right)-2\n_cV_b\left(2\n_c\n_bH_f+2V_a\n_ch_{ab}+h_{ab}\n_cV_a\right)\\
	&\qquad -H^2_f\D^3f(N,N,N)+2H_f\D^3f(\n_aF,N,N)-\D^3f(N,\n_aF,\n_bF)V_aV_b\\
	=&\left(|h|^2+\D^2f(N,N)\right)Z+2X_aU_a+2H_fh_{ab}\left(2\n_th_{ab}-H_fh^2_{ab}+2V_c\n_ch_{ab}\right)\\
	&\qquad -2\n_cV_b\left(2\n_th_{ab}-2H_fh^2_{ab}+2V_a\n_ah_{cb}+h_{ab}\n_cV_a\right)-\D^3f(N,V-H_fN,V-H_fN)\\
	=&\left(|h|^2+\D^2f(N,N)\right)Z+2X_aU_a+2H_f(2W_{ab}-H_fh^2_{ab})-2\n_cV_b(2W_{cb}-2H^2_fh^2_{cb}+h_{ab}\n_cV_a)\\
	&\qquad-\D^3f(N,V-H_fN,V-H_fN)\\
	=&\left(|h|^2+\D^2f(N,N)\right)Z+2X_aU_a-4Y_{ab}W_{ab}+2H_fh^2_{ab}Y_{ab}-2\n_cV_bh_{ab}Y_{ca}\\
	&\qquad -\D^3f(N,V-H_fN,V-H_fN)\\
	=&\left(|h|^2+\D^2f(N,N)\right)Z+2X_aU_a-4Y_{ab}W_{ab}-2Y_{ab}h_{bc}Y_{ca}-\D^3f(N,V-H_fN,V-H_fN).
\end{align*}
Thus we obtain the following proposition:
\begin{proposition}
	The Harnack quantity $Z$ for the $f$-mean curvature flow \eqref{f-mean curvature flow} satisfies the following evolution equation:
	\begin{equation}
		\begin{aligned}
		(\n_t-L)Z=&\left(|h|^2+\D^2f(N,N)\right)Z\\
		&\quad+2X_aU_a-4Y_{ab}W_{ab}-2Y_{ab}h_{bc}Y_{ca}-\D^3f(N,V-H_fN,V-H_fN).
		\end{aligned}
	\end{equation} 
	If $f$ satisfies $\D^3f\equiv 0$, then we have:
	\begin{equation}\label{evolution}
		(\n_t-L)Z= \left(|h|^2+\D^2f(N,N)\right)Z+2X_aU_a-4Y_{ab}W_{ab}-2Y_{ab}h_{bc}Y_{ca}.
	\end{equation}
\end{proposition}

Under the condition $\D^3f\equiv 0$, we prove the Harnack estimate: for all tangent vectors $V$ and $t\geq 0$, $Z(V)\geq 0$.

\begin{theorem}[Differential Harnack Inequality for $f$-Mean Curvature Flow]\label{Harnack estimate theorem}
	Assume $\D^3f\equiv 0$. For any weakly convex solution of the $f$-mean curvature flow \eqref{f-mean curvature flow}, if at $t=0$ for all tangent vectors $V$ we have:
	\begin{equation}\label{Harnack}
		\partial_tH_f+2\langle\n H_f,V\rangle+h(V,V)\geq 0.
	\end{equation}
	Or equivalently, for strictly convex solutions at $t=0$:
	\begin{equation}
	    \partial_t H_f-h^{-1}(\n H_f,\n H_f)\geq 0.
	\end{equation}
	Then for any $t>0$, for all tangent vectors $V$ we have:
	\begin{equation}\label{Harnack}
		\partial_tH_f+2\langle\n H_f,V\rangle+h(V,V)\geq 0.
	\end{equation}
	
\end{theorem}
\begin{proof}
	Given a point $(x_0,t_0)$ and a tangent vector $V$, we can extend $V$ such that at $(x_0,t_0)$ we have $Y=0$ and $U=0$. Then at $(x,t)$, under this extension, we have:
	\[
		(\n_t-L)Z= \left(|h|^2+\D^2f(N,N)\right)Z\stackrel{\triangle}{=}PZ.
	\]
	This shows that $Z$ remains nonnegative. To strictly apply the maximum principle, perturb $Z$ and consider:
	\[
		\overline{Z}=Z+\phi+\psi|V|^2
	\]
	where $\phi, \psi:[0,T]\to\R$ are to be determined.
	
	Assume $\overline{Z}$ is not always positive on $M\times(0,T]$. Then there exists a first time $t_0>0$ and a point $x_0\in M$, $V\in T_{x_0}M$ such that $\overline{Z}=0$. Extend $V$ in a spacetime neighborhood of $(x_0,t_0)$ such that:
	\begin{align*}
		Y_{ab}&=\n_aV_b-H_fh_{ab}=0,\\
		U_a&=(\n_t-L)V_a+h_{ab}\n_bH_f+2\D^2f(\n_aF,N)=0.
	\end{align*} 
	At $(x_0,t_0)$, we have:
	\[
		\partial_t \overline{Z}\leq 0,\quad \n\overline{Z}=0,\quad\Delta\overline{Z}\geq 0.
	\]
	Hence:
	\begin{align*}
		0&\geq(\n_t-L)\overline{Z}= PZ+\partial_t\phi+\partial_t\psi|V|^2+2\psi V_a(\n_t-L)V_a-2\psi|\n V|^2\\
		&\geq P\overline{Z}-P\phi-P\psi|V|^2+\partial_t\phi+\partial_t\psi|V|^2-2\psi h(\n H_f,V)-4\psi \D^2f(V,N)-2\psi H^2_f|h|^2\\
		&\geq (\partial_t-P)\phi+(\partial_t-P)\psi|V|^2-C\psi|V|-C\psi\\
		&\geq (\partial_t-P)\phi+(\partial_t-P-C)\psi|V|^2-C\psi.
	\end{align*}
	Since $|h|^2$ is bounded at time $t_0$, we have $P<C_1$. Taking $\phi=\psi=\epsilon e^{At}$ (with $A>C+C_1$) leads to a contradiction. Therefore, $\hat{Z}\geq 0$ for all $t\geq 0$, and letting $\epsilon\to 0$ yields the theorem.
\end{proof}
\begin{remark}
	When $h>0$, this Harnack inequality is equivalent to:
	\[
		\partial_t H_f-h^{-1}(\n H_f,\n H_f)\geq 0.
	\]
	Its form is consistent with the Li-Yau type Harnack inequality.
\end{remark}

\begin{corollary}
	Under the same assumptions and the additional condition $H_f\geq 0$ at $t=0$, for any solution of the $f$-mean curvature flow \eqref{f-mean curvature flow}, when $t>0$ for all tangent vectors $V$ we have:
	\begin{equation}\label{Harnack1}
		\partial_tH_f+2\langle\n H_f,V\rangle+h(V,V)+\dfrac{H_f}{c(t)}\geq 0.
	\end{equation}
	where $c(t)$ is a positive function for $t\geq 0$.

	In particular, taking $c(t)=2t$ gives the form consistent with Hamilton's mean curvature flow Harnack estimate in \cite{Hamilton95}:
	\begin{equation}\label{Harnack2}
		\partial_tH_f+2\langle\n H_f,V\rangle+h(V,V)+\dfrac{H_f}{2t}\geq 0.
	\end{equation}
\end{corollary}

\section{Integral Harnack Inequality for the $f$-Mean Curvature Flow}

Inequality (\ref{Harnack}) is the differential Harnack inequality for the $f$-mean curvature flow. We can further derive an integral Harnack inequality for $H_f$. Its proof requires the following control conditions.

Let $\lambda,\mu$ be the maximum and minimum eigenvalues of $\n^2_{ab}f$, respectively.
\begin{proposition}\label{Integral inequality}
	Assume $\D^3 f\equiv 0$, $\inf_{M_0}H_f>0$ and $\inf_{M_0}H_f+\lambda-2\mu\leq 0$. Then for any solution of the $f$-mean curvature flow \eqref{f-mean curvature flow}, there exists a constant $C$ depending only on $M_0$ such that for all $t\geq 0$ we have:
	\[
		\dfrac{|h|^2}{H^2_f}\leq C^2.
	\]
\end{proposition}
\begin{proof}
	From the evolution equations:
	\begin{align*}
		\partial_tH_f&=LH_f+(|h|^2+\D^2f(N,N))H_f,\\
		\n_th_{ab}&=Lh_{ab}+(|h|^2+\D^2f(N,N))h_{ab}-h_{ac}\n^2_{cb}f-h_{cb}\n^2_{ac}f.
	\end{align*}
	(where $L=\Delta-\langle\D f,\n\rangle$). Hence:
	\[
		\partial_tH_f=LH_f+(|h|^2+\D^2f(N,N)) H_f\geq LH_f+\mu H_f.
	\]
	Consider the ODE:
	\[
		\frac{dp}{dt}=\mu p,\quad p(0)=\inf_{M_0}H_f.
	\]
	By the ODE comparison principle (see \citep[Corollary 1.4]{Andrews}), we get:
	\[
		\inf_{M_t}H_f(t)\geq p(0)e^{\mu t}\geq \inf_{M_0}H_f,
	\]
	and:
	\begin{align*}
		\partial_t|h|^2&=\n_t|h|^2=2h_{ab}\n_th_{ab}\\
		&=2h_{ab}Lh_{ab}+2(|h|^2+\D^2f(N,N))|h|^2-4h^2_{ab}\n^2_{ab}f\\
		&=L|h|^2-2|\n h|^2+2(|h|^2+\D^2f(N,N))|h|^2-4h^2_{ab}\n^2_{ab}f\\
		&\leq L|h|^2+2|h|^4+(2\lambda-4\mu)|h|^2.
	\end{align*}
	Consider the ODE:
	\[
		\frac{dq}{dt}=2q^2+(2\lambda-4\mu)p,\quad q(0)=\sup_{M_0}|h|^2.
	\]
	Its solution is:
	\[
		q(t)=\frac{(\lambda-2\mu)q(0)e^{(2\lambda-4\mu)t}}{q(0)+(\lambda-2\mu)-2q(0)e^{(2\lambda-4\mu)t}}.
	\]
	Under the condition $q(0)+\lambda-2\mu\leq 0$, the function
	\[
		\theta(x)=\frac{(\lambda-2\mu)q(0)x}{q(0)+(\lambda-2\mu)-q(0)x}
	\]
	is increasing on $[0,1]$, so $q(t)\leq q(0)$ for all $t\geq 0$. By the ODE comparison principle, we get
	\[
		\sup_{M_t} |h|^2(t)\leq p(t)\leq p(0)=\sup_{M_0}|h|^2.
	\]
	Finally:
	\[
		\sup_{M_t}\frac{|h|^2}{H^2_f}(t)\leq \frac{\sup_{M_0}|h|^2}{\inf_{M_0}H^2_f}:=C^2.
	\]
	Since the initial hypersurface $M_0$ is compact or has bounded second fundamental form, the existence of $C$ is obvious.
\end{proof}
\begin{remark}
	The condition $\inf_{M_0}H_f+\lambda-2\mu\leq 0$ is consistent with \cite[Lemma3.3]{JL08}, which also derives the global boundedness of $|h|^2$.
\end{remark}

\begin{corollary}
	Under the conditions of Proposition \ref{Integral inequality}, for any solution of the $f$-mean curvature flow \eqref{f-mean curvature flow} and $t>0$, we have:
	\begin{equation}
		\dfrac{H_f(x_2,t_2)}{H_f(x_1,t_1)}\geq e^{-\frac{C}{4}\Delta},
	\end{equation}
	where $x_1,x_2$ are points on the manifold at times $t_1,t_2$ ($0<t_1<t_2$), and:
	\[
		\Delta=\inf \int\left|\dfrac{dx}{dt}\right|^2_Mdt
	\]
	is the infimum over all paths $x(t)$ satisfying $x(t_1)=x_1, x(t_2)=x_2$. $dx/dt$ is the velocity vector of the path, and $|dx/dt|_M$ is the length of its tangential component on the manifold $M$. In particular:
	\[
		\Delta\leq \frac{d^2(x_1,\hat{x}_2,t_1)}{t_2-t_1},
	\]
	where $d(x_1,\hat{x}_2,t_1)$ is the distance along the manifold from $x_1$ to $\hat{x}_2$ at time $t_1$, and $\hat{x}_2$ is the point that evolves to $x_2$ at time $t_2$.
\end{corollary}

\begin{proof}
	Along the path $y(t)=F(x(t),t)$, we have $H_f(t)=H_f(x(t),t)$, so:
	\[
		\dfrac{d H_f}{dt}=\partial_t H_f+\left\langle\n H_f,\frac{dx}{dt}\right\rangle.
	\]
	Take $V=\frac{1}{2}(dx/dt)^{\top }$ (the tangential component), by the Harnack estimate we have:
	\[
		\dfrac{dH_f}{dt}\geq -\frac{1}{4}h\left(\frac{dx}{dt},\frac{dx}{dt}\right)\geq-\frac{C}{4}H_f\left|\frac{dx}{dt}\right|_M^2.
	\]
	Therefore:
	\[
		\frac{d}{dt}\log H_f\geq-\frac{C}{4}\left|\frac{dx}{dt}\right|_M^2.
	\]
	Integrating yields:
	\[
		\log H_f(x_2,t_2)-\log H_f(x_1,t_1)\geq -\frac{C}{4}\Delta.
	\]
	Taking the exponential of both sides gives the result.
\end{proof}
\begin{remark}
	From (\ref{Harnack1}) we can obtain the integral Harnack inequality:
	\[
		\dfrac{H_f(x_2,t_2)}{H_f(x_1,t_1)}\geq \frac{e^{S(t_1)}}{e^{S(t_2)}}e^{-\frac{C}{4}\Delta}.
	\]
	where $S^{\prime}(t)=1/c(t)>0$. In particular, taking $S(x)=\log(t/2)$ gives the classical form:
	\[
		\dfrac{H_f(x_2,t_2)}{H_f(x_1,t_1)}\geq \sqrt{\frac{t_1}{t_2}}e^{-\frac{C}{4}\Delta}.
	\]
\end{remark}

\bibliographystyle{alpha}
\bibliography{sample}

\end{document}